\documentclass[11pt]{article}
\usepackage{epic,latexsym,amssymb}
\usepackage{tikz}
\usepackage[bold,full]{complexity}
\usepackage{tikz,epsf}

\usepackage{graphicx}[width=\textwidth]

\usepackage{listings}

\usepackage{amsfonts}
\usepackage{amscd}
\usepackage{amsmath}
\usepackage{graphicx}
\usepackage{color}
\usepackage{caption,subcaption}
\usepackage{lineno}

\newtheorem{thm}{Theorem}
\newtheorem{conj}[thm]{Conjecture}
\newtheorem{ob}[thm]{Observation}

\newtheorem{lem}[thm]{Lemma}

\newtheorem{cor}[thm]{Corollary}

\newcommand{\QED}{$\Box$}

\newcommand{\1}{\vspace{0.1cm}}

\newenvironment{proof}{\noindent\textbf{Proof. }}{\hfill $\square$}

\newcommand{\vertex}{\node[vertex]}
\tikzstyle{vertex}=[circle, draw, inner sep=0pt, minimum size=6pt]

\newcommand{\QEDmark}{\mbox{\textsc{qed}}}
\newcommand{\proofStarter}[1]{\textsc{#1} }

\def\vertex(#1){\put(#1){\circle*{2}}}
\def\vertexo(#1){\put(#1){\circle{2}}}
\def\vert(#1){\put(#1){\circle*{1.5}}}
\def\verto(#1){\put(#1){\circle{1.5}}}
\def\lab(#1)#2{\put(#1){\makebox(0,0)[c]{#2}}}
\definecolor{DarkGreen}{rgb}{0.2, 0.6, 0.3}

\definecolor{electricindigo}{rgb}{0.44, 0.0, 1.0}

\begin{document}

\title{Conjectures of TxGraffiti: Independence, domination, and matchings}

\author{$^1$Yair Caro, \, $^{2}$Randy Davila, \, $^3$Michael A. Henning, \, and \, $^2$Ryan Pepper
\\ \\
$^1$Department of Mathematics \\
University of Haifa--Oranim\\
Tivon 36006, Israel \\
\small {\tt e-mail: yacaro@kvgeva.org.il}\\
\\
$^2$Department of Mathematics and Statistics \\
University of Houston--Downtown \\
Houston, TX 77002, USA \\
\small {\tt e-mail: davilar@uhd.edu}\\
\small {\tt e-mail: pepperr@uhd.edu}\\
\\
$^3$Department of Mathematics and Applied Mathematics \\
University of Johannesburg \\
Auckland Park, 2006 South Africa\\
\small \tt e-mail: mahenning@uj.ac.za
\\ 
}

\date{}
\maketitle

\begin{abstract}
TxGraffiti is an automated conjecturing program that produces graph theoretic conjectures in the form of conjectured inequalities. This program written and maintained by the second author since 2017 was inspired by the successes of previous automated conjecturing programs including Fajtlowicz's GRAFFITI and DeLaVi\~{n}a's GRAFFITI.pc. In this paper we prove and generalize several conjectures generated by TxGraffiti when it was prompted to conjecture on the \emph{independence number}, the \emph{domination number}, and the \emph{matching number} (and generalizations of each of these graph invariants). Moreover, in several instances we also show the proposed inequalities relating these graph invariants are sharp. 
\end{abstract}

{\small \textbf{Keywords:}  Automated conjecturing; independence number; domination number; matching number}\\
\indent {\small \textbf{AMS subject classification: 05C69}}

\newpage
\section{Introduction}


Over the last decade there has been a surge in the use of artificial intelligence and machine learning across a wide range of disciplines, with many surprising breakthroughs. For example, machine learning clustering algorithms provide effective web search and deep neural networks allow for self-driving cars. Many of these modern advances in machine learning and artificial intelligence are modeled on the basis of designing intelligent machines to solve real world problems. Though applications of artificial intelligence and machine learning are clearly important, it is worth noting that the founder of artificial intelligence, Alan Turing, suggested in 1948 that designing machines to do mathematical research would be a good starting point for this aim~\cite{Turing}. In considering Turing's suggestion, two fundamental questions arise: Can computers be designed to prove theorems? Can computers be designed to pose meaningful conjectures? The latter of which we provide evidence for in this paper.

The task of programming machines to form mathematical conjectures is referred to as \emph{automated conjecturing}, and the first computer program to make conjectures leading to published mathematical research was Fajtlowicz's GRAFFITI~\cite{Graffiti}. This program written in the early-1990's, considers a small collection of mathematical objects (predominantly simple graphs) and invokes heuristics, that by design limit both the number of output statements of the program and the quality of output statements of the program. One notable result inspired by a conjecture of GRAFFITI and related to graph theory is the result due to Favaron, Mah\'{e}o, and Sacl\'{e}~\cite{FaMaSa-91}, which states that $\alpha(G) \ge R(G)$, where $\alpha(G)$ is the \emph{independence number} of the graph $G$ and $R(G)$ is the \emph{residue} of the graph $G$.

As a successor to GRAFFITI, DeLaVi\~{n}a's GRAFFITI.pc~\cite{Graffiti.pc} also invokes heuristics on a collection of graphs in order to produce graph theoretic conjectures. Unlike GRAFFITI, GRAFFITI.pc maintains and computes data on millions of graphs and has been particularly successful in generating conjectures related to domination type graph invariants. GRAFFITI and GRAFFITI.pc are not the only notable automated conjecturing programs, and we would like to mention Lenat’s AM~\cite{Lenat_1, Lenat_2, Lenat_3}, Epstein’s GT~\cite{Epstein_1, Epstein_2}, Colton’s HR~\cite{Colton_1, Colton_2, Colton_3}, Hansen and Caporossi’s AGX~\cite{AGX_1, AGX_2, AGX_3}, and Mélot’s Graphedron~\cite{graphedron_1}.

In this paper we consider conjectures of the automated conjecturing program \emph{TxGraffiti}. This program written in the Python programming language (versions 3.6 and higher) by the second author in 2017, and subsequence versions written in 2019 and 2020, produces conjectures on simple connected graphs. Currently, TxGraffiti stores all connected graphs on 8 or less vertices, all connected cubic graphs on up to 16 vertices, and several other graphs that appear in the literature. For each graph that is stored in the programs database, TxGraffiti also stores pre-computed numerical and boolean values stored for the graph in question. For example, a graph's order, size, domination number, and independence number are stored, together with properties of the graph such as claw-free, triangle-free, or planar; among many other numerical invariants and structural properties. For a video demo of using TxGraffiti, see~\cite{TxGraffiti}. 

The conjectures produced by TxGraffiti are in the form of proposed inequalities relating graph invariants subject to some structural property. We note that TxGraffiti has already produced conjectures that have led to publishable mathematical results, for example, one conjecture of TxGraffiti inspired the surprising result given in~\cite{CaDaPe2020}, which states $\delta(G)\alpha(G) \le \Delta(G)\mu(G)$ for graphs with minimum degree $\delta(G)$, maximum degree $\Delta(G)$, and \emph{matching number} $\mu(G)$. For a notable open conjecture of TxGraffiti, we mention the conjecture which states $Z(G) \le \alpha(G) + 1$, where $G \neq K_4$ is a connected cubic graph and $Z(G)$ is the \emph{zero forcing number} of $G$; a partial solution to this conjecture is given in~\cite{DaHe_zero_forcing}. 

In the following sections of this paper we present, generalize, and prove conjectures of TxGraffiti when the program was asked to focus on the independence number, the edge domination number, and the connected domination number.  In doing so, we give further evidence for the usefulness of TxGraffiti and other automated conjecturing programs in stimulating mathematical research.  

\noindent\textbf{Graph Terminology and Notation.} Throughout this paper all graphs will be considered undirected, simple, and finite. We in general follow the graph theory notation and terminology in~\cite{HeYe-book}. In particular, the order of a graph $G$ with vertex set $V(G)$ and edge set $E(G)$ is $n(G) = |V(G)|$, whereas the size of $G$ is $m(G) = |E(G)|$. Two vertices $v,w \in V(G)$ are neighbors, or adjacent, if $vw \in E(G)$. The open neighborhood of $v\in V(G)$, is the set of neighbors of $v$, denoted $N_G(v)$, whereas the closed neighborhood is $N_G[v] = N_G(v) \cup \{v\}$. The open neighborhood of $S\subseteq V$ is the set of all neighbors of vertices in $S$, denoted $N_G(S)$, whereas the closed neighborhood of $S$ is $N_G[S] = N_G(S) \cup S$. The degree of a vertex $v\in V$ is denoted by $d_G(v) = |N_G(v)|$. The maximum and minimum degree of $G$ is denoted by $\Delta(G)$ and $\delta(G)$, respectively. For $r \ge 1$ an integer, an $r$-\emph{regular graph} is a graph  with all vertices of degree~$r$.

For a set of vertices $S \subseteq V(G)$, the subgraph induced by $S$ is denoted by $G[S]$. If $v \in V(G)$, we denote by $G - v$ the graph obtained from $G$ by deleting $v$. We denote the path, cycle, and complete graph on $n$ vertices by $P_n$, $C_n$, and $K_n$, respectively. A \emph{triangle} in $G$ is a subgraph of $G$ isomorphic to $K_3$, whereas a \emph{diamond} in $G$ is a subgraph of $G$ isomorphic to $K_4$ with one edge removed. A graph $G$ is \emph{$F$-free} if $G$ does not contain $F$ as an induced subgraph. In particular, if $G$ is $F$-free, where $F = K_{1,3}$, then $G$ is \emph{claw-free}. Claw-free graphs are heavily studied and an excellent survey of claw-free graphs has been written by Flandrin, Faudree, and Ryjacek~\cite{claw-free}.

Two edges in a graph $G$ are \emph{independent} if they are not adjacent in $G$; that is, if they have no vertex in common. For a set $M$ of edges of $G$, we let $V(M)$ denote the set of vertices of $G$ that are incident with an edge in $M$. The set $M$ is a \emph{matching} in $G$ if the edges in $M$ are pairwise disjoint, that is, $|V(M)| = 2|M|$. A matching of maximum cardinality is a \emph{maximum matching}, and the cardinality of a maximum matching of $G$ is the \emph{matching number} of $G$, denote by $\mu(G)$.
%
If $M$ is a matching of $G$, a vertex is \emph{$M$-matched} if it is incident with an edge of $M$; otherwise, the vertex is \emph{$M$-unmatched}.

A matching $M$ in $G$ is \emph{maximal} if it is maximal with respect to inclusion, that is, the set $V(G) \setminus V(M)$ is independent. The minimum cardinality of a maximal matching in $G$ is the \emph{edge domination number} of $G$, denoted by $\gamma_e(G)$. A maximal matching in $G$ of cardinality $\gamma_e(G)$ is a \emph{minimum maximal matching}. Matchings in graphs are extensively studied in the literature (see, for example the classic book on matching by Lovász and Plummer~\cite{LoPl-86}, and the excellent survey articles by Plummer~\cite{Pl-03} and Pulleyblank~\cite{Pu-95}). Edge domination in graphs is well studied with over 125 papers listed on MathSciNet to date. For a recent paper, we refer the reader to~\cite{BaFuHeMoRa-20}.

The \emph{independence number}, $\alpha(G)$, of $G$ is the maximum cardinality of an independent set of vertices in $G$. For $k \ge 0$ an integer, a set $S \subseteq V(G)$ is a $k$-\emph{independent set} in $G$ if $\Delta(G[S]) \le k$. The cardinality of a maximum $k$-independent set in $G$ is the $k$-\emph{independence number} of $G$, denoted by $\alpha_k(G)$. We note that the $0$-independence number of $G$ is the classic independence number, that is, $\alpha(G) = \alpha_0(G)$. The independence and $k$-independence numbers are very well-studied in the literature, and a survey article is given by Chellali, Favaron, Hansberg, and Volkmann~\cite{ChFaHaVo-12}. A refer to~\cite{k-alpha:1, Amos, k-alpha:2, k-alpha:3, k-alpha:4,Fa-84,HeLoSoYe-14} for a few selected papers on $k$-independence.

A set $S$ of vertices in a graph $G$ is a \emph{dominating set} if every vertex not in $S$ is adjacent to a vertex in~$S$. A dominating set $S$ with the additional property that every vertex in $S$ is adjacent to some other vertex in $S$ is a \emph{total dominating set}, abbreviated TD-set of $G$. Moreover, a dominating set that is independent is an \emph{independent dominating set}. A dominating set $S$ of $G$ such that $G[S]$ is connected is a \emph{connected dominating set}. For $k \ge 1$, a $k$-\emph{dominating set} of $G$ is a dominating set $S$ such that every vertex outside $S$ has at least $k$ neighbors in $S$, that is, $|N_G(v) \cap S| \ge k$ for every vertex $v \in V(G) \setminus S$. The \emph{domination number} $\gamma(G)$ of $G$ is the minimum cardinality of a dominating set of $G$, and a dominating set of $G$ of cardinality $\gamma(G)$ is called a $\gamma$-\emph{set of $G$}. The \emph{total domination number}, \emph{independent domination number}, \emph{connected domination number}, \emph{$k$-domination number} are defined analogously, and denoted by $\gamma_t(G)$, $i(G)$, $\gamma_c(G)$ and $\gamma_k(G)$, respectively. For recent books on domination and total domination in graphs, we refer the reader to~\cite{HaHeHe-20,HaHeHe-21,HeYe-book}.

For sets $X$ and $Y$ of vertices in a graph $G$, we denote by $[X,Y]$ the set of edges of $G$ with one end in $X$ and the other end in $Y$. For $\ell \ge 1$ an integer, we use the standard notation $[\ell] = \{1,\ldots,\ell\}$.

\section{Independence versus domination}

In this section we prove and generalize several conjectures of TxGraffiti relating independence and domination in graphs. For similar results relating these two invariants, see also \cite{Caro_Pep, Hansberg_Pepper, Pep}. The first TxGraffiti conjecture we consider relates the independence number to the total domination number.
\begin{conj}[TxGraffiti, confirmed.]

\label{conj:1}
If $G$ is a connected cubic graph, then $\alpha(G) \le \frac{3}{2}\gamma_t(G)$.
\end{conj}

As with many conjectures of TxGraffiti, natural generalizations can be found with some inspection. This is the case for Conjecture~\ref{conj:1}, and we will prove this generalization. However, before providing this result, we first define a class of graphs that were discovered while considering its statement. For $r \ge 2$, let $\mathcal{G}_r$ be the family of all $r$-regular bipartite graphs $G$ for which there exists a partition $V(G) = A \cup B$, where $A$ induces an $(r-1)$-regular graph and $B$ induces a $1$-regular graph.

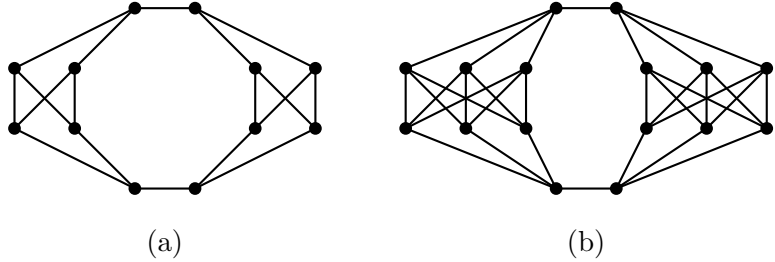
\begin{figure}[htb]
\begin{center}
\begin{tikzpicture}[scale=.8,style=thick,x=1cm,y=1cm]
\def\vr{2.5pt} 
\path (0, 3) coordinate (x1);
\path (1, 3) coordinate (x2);

\path (0, 0) coordinate (y1);
\path (.5, -.5) coordinate (ymid);
\path (1, 0) coordinate (y2);

\path (-2, 2) coordinate (a1);
\path (-1, 2) coordinate (a2);
\path (-2, 1) coordinate (b1);
\path (-1, 1) coordinate (b2);

\path (2, 2) coordinate (c1);
\path (3, 2) coordinate (c2);
\path (2, 1) coordinate (d1);
\path (3, 1) coordinate (d2);

%

\draw (x1) -- (x2);
\draw (y1) -- (y2);

\draw (a1) -- (b1);
\draw (a1) -- (b2);
\draw (a2) -- (b1);
\draw (a2) -- (b2);

\draw (a1) -- (x1);
\draw (a2) -- (x1);
\draw (b1) -- (y1);
\draw (b2) -- (y1);

\draw (c1) -- (d1);
\draw (c1) -- (d2);
\draw (c2) -- (d1);
\draw (c2) -- (d2);

\draw (c1) -- (x2);
\draw (c2) -- (x2);
\draw (d1) -- (y2);
\draw (d2) -- (y2);

\draw (x1) [fill=black] circle (\vr);
\draw (x2) [fill=black] circle (\vr);
\draw (y1) [fill=black] circle (\vr);
\draw (y2) [fill=black] circle (\vr);

\draw (a1) [fill=black] circle (\vr);
\draw (a2) [fill=black] circle (\vr);
\draw (b1) [fill=black] circle (\vr);
\draw (b2) [fill=black] circle (\vr);

\draw (c1) [fill=black] circle (\vr);
\draw (c2) [fill=black] circle (\vr);
\draw (d1) [fill=black] circle (\vr);
\draw (d2) [fill=black] circle (\vr);

\draw[anchor = north] (ymid) node {{(a)}};

\path (7, 3) coordinate (x1);
\path (8, 3) coordinate (x2);

\path (7, 0) coordinate (y1);
\path (7.5, -.5) coordinate (ymid);
\path (8, 0) coordinate (y2);

\path (4.5, 2) coordinate (a1);
\path (5.5, 2) coordinate (a2);
\path (6.5, 2) coordinate (a3);
\path (4.5, 1) coordinate (b1);
\path (5.5, 1) coordinate (b2);
\path (6.5, 1) coordinate (b3);

\path (8.5, 2) coordinate (c1);
\path (9.5, 2) coordinate (c2);
\path (10.5, 2) coordinate (c3);
\path (8.5, 1) coordinate (d1);
\path (9.5, 1) coordinate (d2);
\path (10.5, 1) coordinate (d3);

%

\draw (x1) -- (x2);
\draw (y1) -- (y2);

\draw (a1) -- (b1);
\draw (a1) -- (b2);
\draw (a1) -- (b3);
\draw (a2) -- (b1);
\draw (a2) -- (b2);
\draw (a2) -- (b3);
\draw (a3) -- (b1);
\draw (a3) -- (b2);
\draw (a3) -- (b3);

\draw (a1) -- (x1);
\draw (a2) -- (x1);
\draw (a3) -- (x1);
\draw (b1) -- (y1);
\draw (b2) -- (y1);
\draw (b3) -- (y1);

\draw (c1) -- (d1);
\draw (c1) -- (d2);
\draw (c1) -- (d3);
\draw (c2) -- (d1);
\draw (c2) -- (d2);
\draw (c2) -- (d3);
\draw (c3) -- (d1);
\draw (c3) -- (d2);
\draw (c3) -- (d3);

\draw (c1) -- (x2);
\draw (c2) -- (x2);
\draw (c3) -- (x2);
\draw (d1) -- (y2);
\draw (d2) -- (y2);
\draw (d3) -- (y2);

\draw (x1) [fill=black] circle (\vr);
\draw (x2) [fill=black] circle (\vr);
\draw (y1) [fill=black] circle (\vr);
\draw (y2) [fill=black] circle (\vr);

\draw (a1) [fill=black] circle (\vr);
\draw (a2) [fill=black] circle (\vr);
\draw (b1) [fill=black] circle (\vr);
\draw (b2) [fill=black] circle (\vr);
\draw (a3) [fill=black] circle (\vr);
\draw (b3) [fill=black] circle (\vr);

\draw (c1) [fill=black] circle (\vr);
\draw (c2) [fill=black] circle (\vr);
\draw (d1) [fill=black] circle (\vr);
\draw (d2) [fill=black] circle (\vr);
\draw (c3) [fill=black] circle (\vr);
\draw (d3) [fill=black] circle (\vr);

\draw[anchor = north] (ymid) node {{(b)}};





\end{tikzpicture}
\end{center}
\vskip -0.5 cm
\caption{Graphs in the families $\mathcal{G}_3$ and $\mathcal{G}_4$ }
\label{fig:sharp0}
\end{figure}

For example, the graph $G$ shown in Figure~\ref{fig:sharp0}(a) is a $3$-regular graph that belongs to the family $\mathcal{G}_3$, and the graph $G$ shown in Figure~\ref{fig:sharp0}(b) is a $4$-regular graph that belongs to the family $\mathcal{G}_4$. Further, we note that the graph $G$ in Figure~\ref{fig:sharp0}(a) satisfies $\alpha(G) = \frac{3}{2}\gamma_t(G)$, while the graph $G$ shown in Figure~\ref{fig:sharp0}(b) satisfies $\alpha(G) =  \frac{4}{2}\gamma_t(G) = \frac{1}{2}\gamma_t(G)$.

Suppose that $G \in \mathcal{G}_r$ for some $r \ge 2$. Thus, $G$ is a $r$-regular bipartite graph that contains a partition $V(G) = A \cup B$, where $A$ induces an $(r-1)$-regular subgraph and $B$ induces a $1$-regular graph. Let $G$ have order~$n$, and so $n = |A| + |B|$. Since every vertex in $A \cup B$ has a neighbor in $B$, the set $B$ is a TD-set of $G$, and so $\gamma_t(G) \le |B|$. Doubling counting edges between $A$ and $B$, we have $(r-1)|B| = |A| = n - |B|$, and so $\gamma_t(G) \le |B| = n/r$. Every graph $G$ of order~$n$ without isolated vertices satisfies $\gamma_t(G) \ge n/\Delta(G)$. In our case the graph $G$ is $r$-regular, and so $\gamma_t(G) \ge n/r$. Consequently, $\gamma_t(G) = n/r$. By definition the graph $G \in \mathcal{G}_r$ is a regular bipartite graph, implying that $\alpha(G) = n/2$. As observed earlier, $\gamma_t(G) = n/r$. Consequently, $\alpha(G) = \frac{r}{2}\gamma_t(G)$. We state this formally as follows.

\begin{ob}
\label{ob:1} 
For $r \ge 2$, if $G \in \mathcal{G}_r$, then $\alpha(G) = \frac{r}{2}\gamma_t(G)$. 
\end{ob}

For $r \ge 2$, the independence number and total domination number of an $r$-regular graph are related as follows.

\begin{thm}
\label{thm:Tdom-vs-indep}
For $r \ge 2$, if $G$ is an $r$-regular graph, then
\[
\alpha(G) \le \frac{r}{2}\gamma_t(G),
\]
with equality if and only if $G \in \mathcal{G}_r$.
\end{thm}
\proof For $r \ge 2$, let $G$ be an $r$-regular graph of order~$n$. Let $D$ be a $\gamma_t$-set of $G$, and so $D$ is a TD-set of $G$ and $|D| = \gamma_t(G)$. Further, let $Q = V(G) \setminus D$. Since $D$ is a TD-set of $G$, each vertex in $D$ is adjacent with at most $r-1$ vertices in $Q$. Furthermore since $D$ is a dominating set, every vertex in $Q$ has at least one neighbor in $D$. Thus, the number of edges between $D$ and $Q$ is at most $(r-1)|D|$ and at least $|Q|=n-|D|$, from which we deduce; 
\begin{equation}
\label{Eq:1A}
\gamma_t(G) = |D| \ge \frac{n}{r}.
\end{equation}

As first observed in 1964 by Rosenfeld~\cite{Ro-64}, the independence number of a regular graph is at most one-half its order, that is, $\alpha(G) \le \frac{1}{2}n$. Hence, by Equation~(\ref{Eq:1A}), we have 
\begin{equation}
\label{Eq:1B}
\gamma_t(G) \ge \frac{n}{r} \ge \frac{2}{r}\alpha(G), 
\end{equation}

\noindent 
or, equivalently, $\alpha(G) \le \frac{r}{2}\gamma_t(G)$. This establishes the desired upper bound in the statement of the theorem.

Suppose, next, that $G$ is an $r$-regular graph that satisfies $\alpha(G) = \frac{r}{2}\gamma_t(G)$ where $r \ge 2$. This implies that the inequalities in Inequality Chains~(\ref{Eq:1A}) and~(\ref{Eq:1B}) are all equalities. Therefore, $\gamma_t(G) = |D| = \frac{n}{r}$ and $\alpha(G) = \frac{1}{2}n$. Moreover, every vertex in $D$ is adjacent with exactly $r-1$ vertices in $Q$, and every vertex in $Q$ has exactly one neighbor in $D$ (for otherwise, we would have strict inequality in Inequality~(\ref{Eq:1A})). Hence, the subgraph $G[Q]$ of $G$ induced by the set $Q$ is an $(r-1)$-regular graph, and the subgraph $G[D]$ of $G$ induced by the set $D$ is a $1$-regular graph. Therefore, $|D| =  n/r$ and $|Q| = n(r-1)/r$. Since $\alpha(G) = \frac{1}{2}n$, the graph $G$ is a bipartite graph. Hence, $G$ is an $r$-regular bipartite graph for which there exists a partition $V(G) = A \cup B$, where $A = Q$ induces an $(r-1)$-regular graph and $B = D$ induces a $1$-regular graph. Therefore by definition, we have $G \in \mathcal{G}_r$. Conversely, if $G \in \mathcal{G}_r$, then by Observation~\ref{ob:1}, we have $\alpha(G) = \frac{r}{2}\gamma_t(G)$.~\QED

\medskip 
TxGraffiti also produces conjectures on graphs with forbidden subgraphs. We remark that if a TxGraffiti  presents to the user a conjecture that requires a forbidden subgraph, then the proposed inequality must be false for at least one graph with the forbidden subgraph present. That is, TxGraffiti will only produce an inequality relating graph invariants for the largest possible family of graphs for which the inequality is true. For example, the following conjecture states that the independence number of connected and claw-free graphs is bounded from above by the $2$-domination number, and so, among all connected graphs, TxGraffiti found at least one graph with $K_{1,3}$ as an induced subgraph which does not satisfy this relationship. 
\begin{conj}[TxGraffiti, confirmed.]\label{conj:2}
If $G$ is a connected claw-free graph, then $\alpha(G) \le \gamma_{2}(G)$.
\end{conj}

As with Conjecture~\ref{conj:1}, a natural generalization of Conjecture~\ref{conj:2} was not difficult to deduce from the original proposed inequality, and we prove this generalization with the following theorem.
\begin{thm}\label{thm:main2}
For $r \ge 3$ and $k \ge 0$, if $G$ is a $K_{1,r}$-free graph and $j = r(k+1)-1$, then
\[
\alpha_k(G)\le \gamma_{j}(G),
\]
and this bound is sharp.
\end{thm}
\proof For integers $r \ge 3$ and $k \ge 0$, let $G$ be a $K_{1,r}$-free graph and let $j = r(k+1)-1$. Let $I$ be a maximum $k$-independent set of $G$, and let $D$ be a minimum $j$-dominating set of $G$. That is, $|I| = \alpha_k(G)$ and $|D| = \gamma_j(G)$. Next let $A = I \setminus D$ and $B = D \setminus I$. We note that $I = A \cup (I \cap D)$ and $D = B \cup (I \cap D)$. If $A =\emptyset$, then $I \subseteq D$, and so $\alpha_k(G) = |I| \le |D| = \gamma_j(G)$. Hence, we may assume $A \ne \emptyset$, for otherwise the desired result, namely $\alpha_k(G)\le \gamma_{j}(G)$, is immediate.

Suppose that $B = \emptyset$. If $D \subset I$ and $v$ is an arbitrary vertex in $I \setminus D$, then $v$ has at least~$j$ neighbors in $D$. However since $j = r(k-1) - 1 > k$, this would imply that $v$ has strictly greater than~$k$ neighbors in $D$, contradicting the fact that $I$ is a $k$-independent set of $G$. Hence, in this case when $B = \emptyset$, we must have $D = I$, and so $\alpha_k(G) = |I| = |D| = \gamma_j(G)$. Hence, we may further assume $B \ne \emptyset$, for otherwise the desired result is immediate.

Recall that $D = B \cup (I \cap D)$ and that $A = I \setminus D$. Since $I$ is a $k$-independent set of $G$, every vertex in $A$ is adjacent to at most $k$ vertices in $I \cap D$. Since $D$ is a $j$-dominating set of $G$, every vertex in $A$ is adjacent to at least~$j = r(k+1)-1$ vertices in $D$, implying that every vertex in $A$ is adjacent to at least~$r(k+1)-1-k = (r-1)(k+1)$ vertices in $B$. Thus, the number of edges, $|[A,B]|$, between $A$ and $B$ is given by
\begin{equation}
\label{Eq1}
|[A,B]| \ge (r-1)(k+1)|A|.
\end{equation}

Suppose that some vertex $v \in B$ has at least $(r-1)(k +1) + 1$ neighbors in $A$. Let $A_v = N(v) \cap A$, and let $G_v = G[A_v]$ be the subgraph of $G$ induced by the set $A_v$. Thus, $G_v$ has order~$|A_v| \ge (r-1)(k +1) + 1$. We note that the maximum degree $\Delta(G_v)$ in $G_v$ is at most~$k$. Therefore,
\[
\alpha(G_v) \ge \frac{|V(G_v)|}{\Delta(G_v) + 1} \ge \frac{(r-1)(k+1) + 1}{k+1} > r-1, \1
\]
implying that $\alpha(G_v) \ge r$. Let $I_v$ be a maximum independent set in $G_v$, and so $|I_v| = \alpha(G_v) \ge r$. We note that the graph $G[I_v \cup \{v\}] = K_{1,r}$, contradicting our supposition that the graph $G$ is $K_{1,r}$-free. Hence, every vertex in $B$ has at most $(r-1)(k+1)$ neighbors in $A$, implying that
\begin{equation}
\label{Eq2}
|[A,B]| \le (r-1)(k+1)|B|.
\end{equation}

By Inequalities~(\ref{Eq1}) and~(\ref{Eq2}), we have $|A| \le |B|$. Thus,
\[
\alpha_k(G) = |I| = |A| + |I \cap D| \le |B| + |I \cap D| = |D| = \gamma_j(G).
\]
This completes the proof of Theorem~\ref{thm:main2}. 

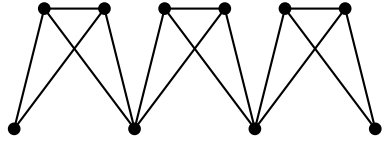
\begin{figure}[htb]
\begin{center}
\begin{tikzpicture}[scale=.8,style=thick,x=1cm,y=1cm]
\def\vr{2.5pt} 
\path (0, 1) coordinate (p11);
\path (1, 1) coordinate (p12);

\path (2, 1) coordinate (p21);
\path (3, 1) coordinate (p22);

\path (4, 1) coordinate (p31);
\path (5, 1) coordinate (p32);

\path (-.5, -1) coordinate (v1);
\path (1.5, -1) coordinate (v2);
\path (3.5, -1) coordinate (v3);
\path (5.5, -1) coordinate (v4);

%

\draw (p11) -- (p12);
\draw (p21) -- (p22);
\draw (p31) -- (p32);

\draw (p11) -- (v1);
\draw (p11) -- (v2);
\draw (p12) -- (v1);
\draw (p12) -- (v2);

\draw (p21) -- (v2);
\draw (p21) -- (v3);
\draw (p22) -- (v2);
\draw (p22) -- (v3);

\draw (p31) -- (v3);
\draw (p31) -- (v4);
\draw (p32) -- (v3);
\draw (p32) -- (v4);

\draw (p11) [fill=black] circle (\vr);
\draw (p12) [fill=black] circle (\vr);
\draw (p21) [fill=black] circle (\vr);
\draw (p22) [fill=black] circle (\vr);
\draw (p31) [fill=black] circle (\vr);
\draw (p32) [fill=black] circle (\vr);
\draw (v1) [fill=black] circle (\vr);
\draw (v2) [fill=black] circle (\vr);
\draw (v3) [fill=black] circle (\vr);
\draw (v4) [fill=black] circle (\vr);




\end{tikzpicture}
\end{center}
\vskip -0.5 cm
\caption{A claw-free graph $G$ with $\alpha(G) = \gamma_2(G)$. }
\label{fig:sharp}
\end{figure}

To see this bound is sharp, consider the family of graphs constructed after the proceeding corollary, one of which is illustrated in Figure~\ref{fig:sharp}.~\QED

\medskip
As a special case of Theorem~\ref{thm:main2}, when $r = 3$ and $k = 0$, we confirm Conjecture~\ref{conj:2} in the affirmative. This result is stated formally with the following corollary, where we further establish sharpness with the infinite family of graphs proceeding its statement. 
\begin{cor}
\label{cor:conj2}
If $G$ is a connected claw-free graph, then 
\[
\alpha(G) \le \gamma_{2}(G),
\] 
and this bound is sharp.
\end{cor}

To see that the bound in Corollary~\ref{cor:conj2} is sharp, let $\ell \ge 1$ and let $G_{\ell}$ be the graph constructed as follows. Take $\ell$ vertex disjoint copies of $P_2$, say $P^1_2, \dots, P^\ell_2$, and $\ell+1$ isolated vertices, say $v_1, v_2, \dots, v_{\ell+1}$, and join both vertices of $P^i_2$ to the two vertices $v_i$ and $v_{i+1}$ for all $i \in [\ell]$. The resulting graph $G_\ell$ is a claw-free graph.
The set $\{v_1, v_2, \dots, v_{\ell+1}\}$ is both a maximum independent set in $G_\ell$ and a $2$-dominating set of $G_\ell$. Therefore, by Corollary~\ref{cor:conj2}, we have $\ell + 1 \le \alpha(G_\ell) \le \gamma_2(G_\ell) \le \ell + 1$. Hence, we have equality throughout this inequality chain, implying that $\alpha(G_\ell) = \gamma_2(G_\ell)$. See Figure~\ref{fig:sharp} for one such example.

For $k \ge 0$, the \emph{local $k$-independence number} $\alpha_{k,L}(G)$ of a graph $G$ is the maximum $k$-independence number among all subgraphs of $G$ induced by the open neighborhoods of the vertices, that is,
\[
\alpha_{k,L}(G) = \max_{v \in V(G)} \left\{ \alpha_k( G[N(v)] ) \right\}.
\]

A second generalization of Conjecture \ref{conj:2}, using a similar proof strategy, and involving local $k$-independence number is presented below. Conjecture \ref{conj:2} is the special case of the theorem below when $r=3$ and $k=0$.

\begin{thm}\label{thm:main3}
For $k \ge 0$ and $r \ge k+2$, if $G$ is a connected graph with $\alpha_{k,L}(G) \le r-1$, then
\[
\alpha_k(G) \le \gamma_{r+k-1}(G).
\]
\end{thm}
\proof For integers $k \ge 0$ and $r \ge k+2$, let $G$ be a connected graph with $\alpha_{k,L}(G) \le r-1$. Let $I$ be a maximum $k$-independent set of $G$, and let $D$ be a minimum $(r+k-1)$-dominating set of $G$. That is, $|I| = \alpha_k(G)$ and $|D| = \gamma_{r+k-1}(G)$. Let $A = I \setminus D$ and $B = D \setminus I$.
If $A =\emptyset$, then $I \subseteq D$, and so $\alpha_k(G) = |I| \le |D| = \gamma_{r+k-1}(G)$. Hence, we may assume that $A \ne \emptyset$. Suppose that $B = \emptyset$. If $D \subset I$ and $v$ is an arbitrary vertex in $I \setminus D$, then $v$ has at least~$r+k-1 \ge k+1$ neighbors in $D$, contradicting the fact that $I$ is a $k$-independent set of $G$. Hence, $D = I$, and so $\alpha_k(G) = |I| = |D| = \gamma_{r+k-1}(G)$. Hence, we may further assume $B \ne \emptyset$, for otherwise the desired result is immediate.

Since $I$ is a $k$-independent set of $G$, every vertex in $A$ is adjacent to at most $k$ vertices in $I \cap D$. Since $D$ is a $(r+k-1)$-dominating set of $G$, every vertex in $A$ is therefore adjacent to at least~$r-1$ vertices in $B$, and so
\begin{equation}
\label{Eq3}
|[A,B]| \ge (r-1)|A|.
\end{equation}

Suppose that some vertex $v \in B$ has at least $r$ neighbors in $A$. Since the local $k$-independence number of $G$ is at most~$r-1$, this implies that the set of neighbors of $v$ in $A$ induces a subgraph of $G$ with maximum degree strictly greater than~$k$. However, such a subgraph of $G$ is a subgraph of the graph $G[I]$ which has maximum degree at most~$k$, a contradiction. Hence, every vertex in $B$ has at most $(r-1)$ neighbors in $A$, implying that
\begin{equation}
\label{Eq4}
|[A,B]| \le (r-1)|B|.
\end{equation}

By Inequalities~(\ref{Eq3}) and~(\ref{Eq4}), we have $|A| \le |B|$. Thus,
\[
\alpha_k(G) = |I| = |A| + |I \cap D| \le |B| + |I \cap D| = |D| = \gamma_{r+k-1}(G).
\]
This completes the proof of Theorem~\ref{thm:main3}.~\QED

While both Theorems \ref{thm:main2} and \ref{thm:main3} are generalizations of Conjecture \ref{conj:2}, and happen to coincide in the case that $k=0$, they become different for higher values of $k$. In general, Theorem \ref{thm:main3} has a better bound on $\alpha_k(G)$ than Theorem \ref{thm:main2}, at the cost of a stronger hypothesis. This means that by excluding more subgraphs than just $K_{1,r}$, Theorem \ref{thm:main3} yields a better upper bound than Theorem \ref{thm:main2}. For example, when $k=3$ and $r=5$, Theorem \ref{thm:main2} gives the result that $\alpha_3(G) \leq \gamma_{19}(G)$, provided $G$ is $K_{1,5}$-free. In the same conditions on $k$ and $r$, Theorem \ref{thm:main3} gives the result that $\alpha_3(G) \leq \gamma_7(G)$, provided $G$ has maximum local $3$-independence number of at most 4 (which means that $G$ is $K_{1,5}$-free, in addition to having many other forbidden subgraphs). 


\section{Edge domination versus matching}

In this section we prove and generalize a conjecture of TxGraffiti on the edge domination number of a graph. Recall that the edge domination number $\gamma_e(G)$ of a graph $G$ is the minimum cardinality of maximal matching in $G$, and a minimum maximal matching is a maximal matching in $G$ of cardinality $\gamma_e(G)$. The following conjecture of TxGraffiti inspired this sections main result.
\begin{conj}[TxGraffiti, confirmed.]\label{conj:3}
If $G$ is a connected cubic graph, then $\gamma_e(G) \ge \frac{3}{5}\mu(G)$.
\end{conj}

While investigating the bound posed in Conjecture~\ref{conj:3}, we subsequently established a more general and stronger lower bound on the edge domination number in terms of the order, maximum degree, and minimum degree. This result is given formally by the following theorem.

\begin{thm}\label{thm:minmax_matching}
If $G$ is a graph of order~$n$, minimum degree~$\delta \ge 1$, and maximum degree $\Delta$, then
\[
\gamma_e(G) \ge \frac{\delta n}{2(\Delta + \delta -1)}.
\]
\end{thm}
\proof Let $G$ be a graph a graph of order~$n$, minimum degree~$\delta \ge 1$, maximum degree $\Delta$, and let $M$ be a minimum maximal matching in $G$. Let $A$ be the set of $M$-matched vertices, and so, $A = V(M)$ and $|A| = 2|M| = 2\gamma_e(G)$. Let $B = V(G) \setminus A$, and so $|B| = n - 2\gamma_e(G)$. By the maximality of the matching $M$, the set $B$ is an independent set.
By double counting the number of edges, $|[A,B]|$, between $A$ and $B$ we have
\begin{equation*}
\label{Eq5}
\delta( n - 2\gamma_e(G) ) = \delta |B| \le |[A,B]| \le (\Delta - 1)|A| = 2(\Delta - 1)\gamma_e(G),
\end{equation*}
or, equivalently, $\gamma_e(G) \ge \delta n/(2(\Delta + \delta -1))$.~\QED

\medskip
If $G$ is a graph of order~$n$, then $\mu(G) \le \frac{1}{2}n$. A \emph{perfect matching} is a matching that matches every vertex; that is, if a graph $G$ of order~$n$ has a perfect matching, then $n$ is even and $\mu(G) = \frac{1}{2}n$. As an immediate consequence of Theorem~\ref{thm:minmax_matching} we have the following results.

\begin{cor}\label{cor:edge-dom-1}
If $G$ is a graph of order~$n$, minimum degree~$\delta \ge 1$, and maximum degree $\Delta$, then
\[
\gamma_e(G) \ge \left( \frac{\delta}{\Delta + \delta -1} \right) \mu(G).
\]
\end{cor}

\begin{cor}\label{cor:edge-dom-2}
For $r \ge 1$, if $G$ is an $r$-regular graph, then $\gamma_e(G) \ge \left( \frac{r}{2r -1} \right) \mu(G)$.
\end{cor}

In the special case when $r=3$ in the statement of Corollary~\ref{cor:edge-dom-2} we have the following result, which resolves Conjecture~\ref{conj:3} in the affirmative.

\begin{cor}\label{cor:edge-dom-3}
If $G$ is a connected cubic graph, then $\gamma_e(G) \ge \frac{3}{5}\mu(G)$.
\end{cor}

\section{Connected domination versus matching and independence}

Next we investigate the connected domination number. Our main results concerning connected domination were inspired by the following two TxGraffiti conjectures, both of which we resolve in the affirmative as corollaries to more general results we present in this section.

\begin{conj}[TxGraffiti, confirmed.]
\label{conj:mu-connected-dom}
If $G$ is a connected cubic graph, then $\gamma_c(G) \ge \mu(G) - 1$.
\end{conj}

\begin{conj}[TxGraffiti, confirmed.]
\label{conj:alpha-connected-dom}
If $G$ is a connected cubic graph, then $\gamma_c(G) \ge \alpha(G) - 1$.
\end{conj}

Conjecture~\ref{conj:alpha-connected-dom} inspired the following more general theorem focusing on the connected domination number. We remark that the following result first appeared in~\cite{CaDaPe2015} as a consequence of an upper bound on the \emph{zero forcing number} of a graph. The proof we provide is different in that our proof technique is focused only on the connected domination number. Moreover, as a consequence of our proceeding arguments we are also able to establish a characterization of equality in the following bound.
\begin{thm}\label{thm:connected-dom}
For $r\ge 2$, if $G$ is a connected $r$-regular graph of order~$n$, then
\[
\gamma_c(G) \ge \frac{n-2}{r-1}.
\]
\end{thm}

\proof For $r \ge 2$, let $G$ be a connected $r$-regular graph with order~$n$. If $G$ contains a dominating vertex, then $G = K_{r+1}$. In this case, $n = r+1$ and $\gamma_c(G) = 1 = \frac{n-2}{r-1}$. Hence, we may assume that $G$ does not contain a dominating vertex, and so, $\gamma_c(G) \ge 2$. Let $D$ be a $\gamma_c$-set of $G$, and so $D$ is a connected dominating set of $G$ such that $|D| = \gamma_c(G) \ge 2$. Since the subgraph $G[D]$ induced by $D$ is connected, it contains a spanning tree $T$ of size $|D| - 1$. Let $A = V(G) \setminus D$. Each vertex in $A$ has at least one neighbor in $D$. Thus by double counting the number of edges, $|[A,D]|$, between $A$ and $D$, we have by the $r$-regularity of $G$ that
\[
\begin{array}{lcl}
n - |D| = |A| & \le & |[A,D]| \1 \\
& = & \displaystyle{ \sum_{v \in D} d_G(v) - 2|E(G[D])| } \1 \\
& \le & r|D| - 2|E(T)| \1 \\
& = & r|D| - 2(|D|-1), \1
\end{array}
\]
or, equivalently, $(r-1)|D|\ge n-2$, and so $\gamma_c(G) = |D| \ge \frac{n-2}{r-1}$.~\QED

\medskip
For $r \ge 3$, we define a tree $T$ to be a $(1,r)$-tree if every vertex of $T$ has degree $1$ or $r$. Let $L(T)$ denote the set of all leaves in a tree $T$, and let $\ell(T) = |L(T)|$, and so $\ell(T)$ is the number of leaves in the tree $T$. We proceed further with the following lemma. Again, we need remark that the following result first appeared in~\cite{CaDaPe2015} as a consequence of an upper bound on the \emph{zero forcing number} of a graph. However, and as before, our proof technique presented here is solely focused on the maximum number of leaves and excludes any notion of zero forcing in graphs.
\begin{lem}
\label{lem:tree}
For $r \ge 3$, if $T$ is a tree of order~$n$ and maximum degree $\Delta(T) \le r$, then
\[
\ell(T) \le \frac{(r-2)n+2}{r-1},
\]
with equality if and only if $T$ is a $(1,r)$-tree.
\end{lem}
\proof For $r \ge 3$, let $T$ be a tree of order~$n$, size~$m$ and maximum degree $\Delta(T) \le r$. Let $n_i$ denote the number of vertices of degree $i$ in $T$ for $i\in [r]$. Thus, \begin{equation}
\label{Eq6}
\sum_{i=1}^{r}n_i = n \hspace*{0.5cm} \mbox{and} \hspace*{0.5cm} \sum_{i=1}^{r}i \cdot n_i = 2m = 2(n-1),
\end{equation}

and so
\[
\sum_{i=2}^{r}(i-1) n_i = n-2.
\]

Hence,
\[
(r-1)n_r = n-2 -\sum_{i=2}^{r-1}(i-1)n_i \le n-2,
\]

and so,
\begin{equation}
\label{Eq7}
n_r \le \frac{n-2}{r-1}.
\end{equation}

By (\ref{Eq6}), we have
\[
\sum_{i=1}^{r}(i-(r-1)) n_i = 2(n-1)-(r-1) n,
\]
or, equivalently,
\begin{equation}
\label{Eq8}
\sum_{i=1}^{r}(r-i-1) n_i = (r-3) n + 2.
\end{equation}

By Inequality~(\ref{Eq7}) and Equation~(\ref{Eq8}), we have
\begin{equation*}
\begin{split}
(r-2)n_1 & =(r-3)n + 2 + n_r + \sum_{i=2}^{r-1}(i+1-r) n_i \\
      & \le (r-3)n + 2 + \frac{n-2}{r-1} + \sum_{i=2}^{r-1}(i+1-r) n_i \\
      & \le \frac{(r-2)^{2}\cdot n + 2(r-2)}{r-1},
\end{split}
\end{equation*}
noting that $(i+1-r)n_i \le 0$ for all $i \in \{2, \dots, r-1\}$. Thus,
\begin{equation}
\label{Eq9}
\ell(T) = n_1 \le \frac{(r-2)n+2}{r-1}.
\end{equation}

Further, if we have equality in Inequality~(\ref{Eq9}), then we must have equality throughout the above inequality chains, implying that $n_i=0$ for $i\in\{2, \dots, r-1\}$. Hence, if we have equality in Inequality (8), then $T$ is a $(1, r)$-tree.

Conversely, if $T$ is a $(1, r)$-tree, then $n_1+n_r=n$ and $n_1+r \cdot n_r = 2m = 2n-2$, and so
\[
n_r=\frac{n-2}{r-1},
\]
which implies that
\[
\ell(T)=n_1=\frac{(r-2)n+2}{r-1}.
\]
This completes the proof of the lemma.~\QED

\medskip
As a consequence of Lemma~\ref{lem:tree}, we have the following result.

\begin{thm}\label{thm:connected-dom-main}
For $r \ge 3$ if $G$ is a connected graph of order~$n$ and maximum degree $\Delta(G) \le r$,  then
\[
\gamma_c(G) \ge \frac{n-2}{r-1},
\]
with equality if and only if $G$ has a spanning $(1,r)$-tree.
\end{thm}
\proof For $r \ge 3$, let $G$ be a connected graph of order~$n$ and maximum degree $\Delta(G) \le r$. We note that $\gamma_c(G) = n - \ell(T)$, where $T$ is a spanning tree of $G$ with the maximum number of leaves. Since $\Delta(G)\le r$, we note that $\Delta(T)\le r$. Thus, by Lemma~\ref{lem:tree}, we have
\[
\gamma_c(G) = n-\ell(T) \ge n - \frac{(r-2)n+2}{r-1} = \frac{n-2}{r-1},
\]
with equality if and only if $T$ is a spanning $(1,r)$-tree in $G$.~\QED

We remark that recognizing if $G$ has $(1,r)$-spanning tree is NP-Complete for $r = 2$, since this amounts to Hamiltonian path problem as can be seen in~\cite{NPC}. On the other hand, it is easy to construct such regular graphs with spanning $( 1,r)$-tree. To see this construction, begin with any $(1,r)$-tree $T$, having $n$ vertices  and  $n_1 =  \frac{n(r-2) +2}{( r-1)}$ are leaves. Next take any $(r -1)$-regular graph $G$ on $n_1$ vertices, which may not necessarily connected, and embed $G$ in the set formed by the leaves of $T$ to obtain the desired $r$-regular graph with a spanning $(1,r)$-tree.

\medskip
As observed in~\cite{Ro-64}, if $G$ is a regular graph of order $n$ without isolated vertices, then $\alpha(G) \le \frac{1}{2}n$. As observed earlier, if $G$ is a graph of order~$n$, then $\mu(G) \le \frac{1}{2}n$. Hence as a consequence of Theorem~\ref{thm:connected-dom-main}, we have a lower bound on the connected domination number of a connected regular graph in terms of the matching number and regularity. Furthermore, Caro et al.~\cite{CaDaPe2020} recently proved $\mu(G)\ge\alpha(G)$ for all $r$-regular graphs with $r\ge 2$. Thus, we have the following result.
\begin{cor}
\label{cor:connected-dom-mu}
For $r\ge 2$, if $G$ is a connected $r$-regular graph, then
\[
\gamma_c(G) \ge \frac{2\mu(G)-2}{r-1}  \ge  \frac{2\alpha(G)-2}{r-1}.
\]
\end{cor}

In the special case of Corollary~\ref{cor:connected-dom-mu} when $r = 3$, we have that if $G$ is a connected cubic graph, then $\gamma_c(G) \ge \mu(G) - 1$ and $\gamma_c(G) \ge \alpha(G) - 1$. This proves Conjectures~\ref{conj:mu-connected-dom}  and~\ref{conj:alpha-connected-dom}. We state this result formally as follows.

\begin{cor}
If $G$ is a connected cubic graph, then the following holds. \1 \\
\indent {\rm (a)} $\gamma_c(G) \ge \mu(G) - 1$. \\
\indent {\rm (b)} $\gamma_c(G) \ge \alpha(G) - 1$.
\end{cor}

\end{document}